\documentclass{elsarticle}

\usepackage{amsmath,amssymb, amsthm}
\usepackage[utf8]{inputenc}
\usepackage{mathrsfs}
\usepackage[margin=2.2cm]{geometry}

\usepackage{graphicx}
\usepackage{subcaption}

\usepackage{wrapfig}
\usepackage{float}

\usepackage{verbatim}

\usepackage{hyperref}

\usepackage{xcolor}
\hypersetup{
	colorlinks   = true, %Colours links instead of ugly boxes
	urlcolor     = {blue!90!black}, %Colour for external hyperlinks
	linkcolor    = {blue!90!black}, %Colour of internal links
	citecolor   = {red!90!black} %Colour of citations
}

\newcommand{\R}{\mathbb{R}}
\newcommand{\N}{\mathbb{N}}
\renewcommand{\d}{\mathrm{d}} 

\theoremstyle{plain}
\newtheorem{thm}{Theorem}[section]
\newtheorem{prop}[thm]{Proposition}
\newtheorem{lem}[thm]{Lemma}

\newtheorem{manualLeminner}{Lemma}
\newenvironment{manualLem}[1]{%
	\renewcommand\themanualLeminner{#1}%
	\manualLeminner
}{\endmanualLeminner}

\theoremstyle{definition}

\newtheorem{con}[thm]{Condition}
\newtheorem{defi}[thm]{Definition}
\newtheorem{rmk}[thm]{Remark}

\newtheorem{example}[thm]{Example}

\begin{document}
	
	\begin{frontmatter}
		
		\title{ { Minimal controllability time for  systems with nonlinear
			drift under a compact convex state constraint \tnoteref{mytitlenote}} }
		\tnotetext[mytitlenote]{ }
		
		% Group authors per affiliation:
		\author[1]{Viktor Bezborodov}
		\author[2]{Luca Di Persio}
		\author[2]{Riccardo Muradore}
		\address[1]{%\emph{Bielefeld University, Faculty of Mathematics} \\
			{Wrocław University of Science and Technology, Faculty of Electronics, Janiszewskiego 11, Wrocław, Poland}}
		\address[2]{Department of Computer Science, University of Verona, Strada le Grazie 15, Verona, Italy}

		\begin{abstract}
			In this paper
			we estimate the minimal controllability time
			for a class of non-linear control systems
			with a bounded convex state  constraint.
			An explicit expression is given for the controllability 
			time 
			 if the 
			image of the control matrix
			is of co-dimension one.
			A lower bound for the controllability time 
			is given in the general case.
			The  technique  is based on finding 	
			a lower dimension system
			with the similar controllability properties as the original system.	
			The controls corresponding to  
			the minimal time, or time close to the minimal one, are discussed and computed analytically.
			The effectiveness of the proposed approach
			is illustrated by a  few examples.	
		\end{abstract}
		
		
	\end{frontmatter}

	%An
	%excerpt 
	%from the Guide for Authors: \emph{
	% Authors are encouraged to prepare their manuscripts as LaTeX documents. For the first submission the use of the standard LaTeX "article" style or the special Automatica style (autart.cls) available on Pampus at
	% www.autsubmit.com is recommended. Authors must not create their own macros or define new commands.
	%}
	
	\section{Introduction}

	In this paper we consider the question of controllability for
	systems with non-linear drift,
	 linear control, and  state constraint. 
	The state of the system is required to stay within certain bounded convex  set.
	
	The proposed
	 technique consists in considering an auxiliary system of lower dimension
	which has similar controllability properties.
	This allows to derive lower bounds on the controllability time.
	For the case when the range of the control matrix has co-dimension one (that is, the image of the control matrix
	is a linear space of dimension one less than the entire state space),
	an explicit expression
	for the controllability time
	is given.
		Using similar technique, in the complementary case
	we give a lower bound on the controllability time. 
	%We also provide three examples. 
	The main idea behind our analysis is 
	that the controllability time 
	for the original system can be expressed
	in terms of the controllability time
	for a lower dimension system.
		The present work is inspired and
		partially
		 motivated by \cite{LTZ18}.

	As in \cite{LTZ18} we focus on controllability with a state constraint,
	but without control constraints, that is, every $L^\infty$ control is allowed.
	Some of the  main techniques in \cite{LTZ18}
	are Brunkovsky normal form for a linear equation,
	and Goh  transformation. In the present paper we too 
	use  equivalent systems to derive properties of
	the 
	minimal controllability time, although our approach differs 
	as the alternative system we arrive to is obtained via orthogonal projection rather than 
	transforming the system into a normal form.
	
	Control systems with state constraints
	is a challenging topic for mathematical analysis that
	has seen a gradual rise in interest over the recent years.
	Quoting from \cite{LTZ18},
	``Controllability under state constraints has not been much investigated
	in the literature, certainly due to the difficulty of the
	question, even for linear control systems.''
	The main object of  \cite{le2017small} is to give conditions on 
	a closed set $S$ so that 
	every point sufficiently close to
	$S$ can be steered into $S$ within a small
	time by an admissible control. 
	The authors call this property small-time local attainability. 
	The control system in \cite{le2017small} is non-linear.
	A similar problem in stochastic settings was studied in \cite{Buckdahn2004}.
	The estimators for systems  with linear and non-linear
	state constraints are surveyed in \cite{simon2010kalman},
	see also \cite{ko2007state}.
	In \cite{krastanov2008constrained},  systems
	with linear state constraint and with convex cone signal constraint 
	are considered. A geometric necessary and sufficient small time controllability 
	condition is formulated in terms of involved constraining sets. 
	Controllability of the fractional systems 
	with constrained delayed controls is treated in \cite{frac_syst1,frac_syst2}.

	We work in a framework similar to \cite{LTZ18}.
	The main differences in the models between \cite{LTZ18}
	and the present paper
	 are that our system is non-linear,
	and that we work only with bounded convex constraint sets.
	In \cite{LTZ18} the focus is on whether the system is controllable under the state constraint,
	and whether the controllability time is positive.
	Meanwhile
	in the present work we
	mostly address the questions
	of estimating and explicitely computing the controllability time
	for non-linear systems in arbitrary dimension.
		It was noted in \cite{LTZ18}
	that obtaining an explicit expression or even an estimate
	for the controllability time
	remained an open problem 
	for linear system in dimensions higher than two.
	Here we provide such expressions
	 and estimates in a wide range of cases encompassing
	 bounded convex state constraints
	 for  more general non-linear systems.
	The main idea behind our analysis is to
	 	show a certain equivalence between the original system and a  non-linear one
	 	with a lower dimension.
	 This is achieved by  decomposing  $\R ^n$
	 into an orthogonal sum  of `fast' directions (those in the range of the control matrix)
 and `slow' directions (the orthogonal complement).
     Under the assumption of convexity the `fast' directions are
  usually straightforward to handle, and the focus of our analysis is on the `slow' ones.
	
	The paper is organized as follows.
	In Section \ref{sec problem form} we describe the non-linear control system
	analyzed in this paper.
	In Section \ref{sec theoretical part} we show 
	that the controllability time of the original system 
	is related to
	the controllability time
	of
	 a lower dimension system.
	 The estimates on the controllability time
	and  an exact expression  are also derived in Section \ref{sec theoretical part}.
	In Section \ref{sec examples} we discuss some numerical examples.
	The concluding remarks  
	 are collected in Section \ref{sec conclusion}.

	\section{Problem formulation} \label{sec problem form}
	
	We consider the system governed by the equation 
	\begin{align}\label{the model}
	\dot{y}(t) & = F( y(t)) + B u (t), 
	\end{align}
	where $y \in \R ^n$ is the state vector,
	$F: \R^n \to \R ^n$ is a continuous vector field, $B$ is $n \times m$ matrix of rank $m$ with $m < n$.
	Henceforth we identify a matrix with the linear operator it induces.
	\begin{comment}
	They satisfy the Kalman condition
	\begin{equation}\label{Kalmal?}
	\text{rank}[B, AB, ..., A ^{n-1}B] = n 
	\end{equation}
	\end{comment}
	\newcommand{\C}{\mathcal{C}}
	System \eqref{the model} is endowed with the additional constraint 
	\[
	y(t )  \in \mathcal{C},
	\]
	where $\mathcal{C} \subset \R ^n$ is a bounded convex  set.
	We always assume  the interior $\C ^{\mathrm{o}}$ 
	of $\C$
	to be not empty.

	Let $y^0, y ^1 \in \C ^{\mathrm{o}}$.
	We define the smallest time to reach one point from another as
	\begin{equation} 
	\begin{split}
	T_\C (y^0, y^1) = \inf\{T >0: \text{there exists }  \label{def T_C}
	u \in L ^\infty([0,T], \R ^m) \text{ s.t. }
	y(t) \in \mathcal{C} \text{ for } t \in [0,T], \text{ and }
	\\
	\eqref{the model} \text{ holds with } 
	y(0)= y^0, y(T) = y^1\}.
	\end{split}
	\end{equation}
	
	Here and throughout, we assume that $y^0, y ^1 \in \C ^{\mathrm{o}}$
	and  adopt the convention $\inf \varnothing = +\infty$.
	If $T_\C (y^0, y^1) < \infty$,
	we say that $y^1$ is \emph{reachable}
	from $y^0$
	with the state constraint $\C$.
	%Since the Kalman condition \eqref{Kalmal?} is satisfied,
	%the controllability time for the system without constraint is always zero for any two points.
	We avoid the initial  and final points being on the boundary ($y^0 \in \partial \C$ or 
		$y^1 \in \partial \C$), because this case would require additional technical assumptions. Indeed, for some
	systems any solution to \eqref{the model} started from some $y^0 \in \partial \C$
   leaves $\C$ immediately. Similarly, for some $y^1 \in \partial \C$
any solution reaching $y^1$ may have to come from the complement of $\C$.
 		On the other hand, if for  $y^0 \in \partial \C$ there exists a signal $u$
 		such that
 		the solution to \eqref{the model} belongs to $\C ^{\mathrm{o}}$
 	 for small $t>0$, then  our results are applicable because 
 	we can take a new starting point  in the interior of $\C$
 	after an arbitrary small delay.

	Let us see in a simple case 
		 a state constraint  affects the controllability time.
	Let $n = 2$, $m = 1$,  $A = \begin{pmatrix}
	0 &  -1 \\ 1 & 0
	\end{pmatrix}$, $B = \begin{pmatrix}
	1 \\ 0
	\end{pmatrix}$, $F(y) = Ay$, so that the system is 
	\begin{equation} \label{example1}
	\begin{cases}
	\begin{aligned}
	\dot{y _1}(t) & = -y_2(t) + u(t),
	\\
	\dot{y _2}(t) & = y_1(t).
	\end{aligned}
	\end{cases}
	\end{equation}
	The Kalman condition is satisfied here, so the state $(0,1)^\top$ (here and elsewhere, $\top$
	indicates transposition)
	can be reached from $(0,0)^\top$  in an arbitrary time if there is no constraint.
	Assume   we also require that for a constant $C >0$,
	$$|y_1(t)| < C.$$
	Then  if system \eqref{example1} reaches $(0,1)^\top$  from $(0,0)^\top$ 
	at time $T>0$, we have
	\[
	1 =y_2(T) - y_2(0) = \int\limits _0 ^ T y_1(t)dt \leq CT,
	\]
	and hence $T \geq \frac 1C$. This means that the controllability time
	cannot be made arbitrary small under the state constraint,
	even though every state within the constraint
	set is reachable 
	from any other.
	More examples can be found in \cite{LTZ18}.

	\begin{comment}
	Another example (given in \cite{LTZ18}) is the system
	\begin{equation}
	\dot{y _1}(t) = y_2 (t),  \ \ \ \dot{y _2}(t) = u (t),
	\end{equation}
	under the constraint $y_2 (t) \geq 0$. Citing from \cite{LTZ18},
	``For any trajectory, the first component $y_1(t)$ must be
	nondecreasing, and then obviously one cannot pass from any point
	to any other."
	\end{comment}
	
			Let 
	$H = (\text{ran}(B))^\perp$ be the orthogonal complement
	to the range of $B$ in $\R^n$.
	The space $H$ represents the `slow' directions mentioned in the introduction.
	Note that $\dim H = n - m$.
	 For a subspace $G $, let $P_{G}y$ be the orthogonal projection of $y$ on $G$ 
	and for a map $M: \R^n \to \R ^n$, $M_{G}(y) := P_{G}M(y)$.
	Denote also by $\C _H $ the orthogonal projection of $\C $ on $H$.
	For a map $M:  X \to Y$ and a set $\mathcal{Q} \subset X$, the image of $\mathcal{Q}$
   under $M$ is defined as
	$
	M \mathcal{Q}= \{ Mx : x \in \mathcal{Q} \} \subset Y.
	$
	
	We make the following assumptions on $F$ and $\C$.
	\begin{con}	\label{con on C}
		The set $\C$ is a bounded  convex
		subset of $\R ^n$ with a smooth $C^1$ boundary.
	\end{con}
	\begin{con} The function $F$ is  continuous and Lipschitz with the Lipschitz constant $L_{{ F }} >0$,
	that is, ${|F(x) - F(y)| \leq L_F|x-y|}$ for all $x, y \in \R^n$. 
	\end{con}
	We also make the following technical assumptions.
	
	\begin{con}[measurable selection] \label{con measurable selection}
			There exists a Borel measurable map
			$f$
		defined on 
		$$\mathcal{D}_f := \{(h_1, h_2) \in H \times H:  \text{ for some } h^\perp \in H^\perp,
		h_1 = F_H(h_2 + h^\perp) \text{ and }
		h_2 + h^\perp \in \C \}$$
		such that for every $(h_1, h_2) \in \mathcal{D}_f$,
		$$
		h_1 = F_H(h_2 + f(h_1, h_2) ). 
		$$
	\end{con}

	Condition \ref{con measurable selection} is a technical assumption which we expect 
	to hold in all reasonable cases. The measurable selection property
	is closely related to the uniformization problem in descriptive set theory \cite{Desc_set_theory}. 
	In particular, if for each $(h_1, h_2) \in \mathcal{D}_f$ the set 
	$$S_{(h_1, h_2)} = \{ h^\perp \mid
			h_1 = F_H(h_2 + h^\perp) \text{ and }
	h_2 + h^\perp \in \C \} $$
    is at most countable or is of positive Lebesgue measure,
    Condition \ref{con measurable selection} is satisfied \cite{surv_meas_selec, meas_selec_large_sec}.
	
	For vectors $v_1,v_2 $ of equal dimension denote by $[v_1, v_2]$ their closed
convex hull. For $x \in \R ^\d$, $\d \in \N$, and   $r>0$, 
let $\mathcal{B}(x,r)$ be the closed ball $\{y \in \R ^\d \mid |y-x|  \leq r \}$.
In particular, $\mathcal{B}(x, 0) = \{x\}$.
	
	\section{Reduction to a lower dimension problem}\label{sec theoretical part}

	One of the aims of the present work is to find another representation of $T_\C (y^0, y^1)$
	as a certain time related to a problem in lower dimension. To this end we introduce 
	 auxiliary dynamics defined by the inclusion
	\begin{align}\label{charade}
	\dot z(t) \in F_{H }\left((z(t) + H^\perp) \cap \C \right),
	\end{align}
	where the state $z$ takes values in $\R ^{\dim (H)}$.
	We now define the controllability time for \eqref{charade} by
	
	\begin{equation} \label{faucet}
	\begin{split}
	\overline T_\C (y^0, y^1) = 
	\inf\big\{T: \text{there exists } 
	z(t) \in C([0,T], H ) 
	\text{ s.t. } \eqref{charade} \text{ holds with } 
	\\
	z(0) = P_{H}y^{0},
	z(T) = P_{H}y^{1}  \big\}.
	\end{split}
	\end{equation}

	Denote by $L ^0(X,Y)$ the set of all measurable maps from $X$ to $Y$.
	We also define another auxiliary  equation with constraints 
	\begin{equation}\label{system hat}
	\begin{gathered}
	\dot z(t) = F_{H }(z(t) + h^\perp(t))
	\\
	z(t) + h^\perp(t) \in \C, \ \ \  h^\perp(t) \in  H ^\perp, 
	\end{gathered}
	\end{equation}
	and the respective controllability time 
	\begin{equation}
	\begin{split}
	\widehat T_\C (y^0, y^1) =  \label{faucet one}
	\inf\big\{T: \text{there exists } z(t) \in C([0,T], H ),
	h^\perp(t) \in L^0([0,T], H ^\perp )
	\text{ s.t. } 
	z(0) = P_{H}y^{0},
	\\
	z(T) = P_{H}y^{1}, 
	\text{ and \eqref{system hat} holds}  \big\}. \notag
	\end{split}
	\end{equation}

	The relation between $ T_\C (y^0, y^1) $, $\overline T_\C (y^0, y^1) $,
	and $\widehat T_\C (y^0, y^1) $ is clarified in this section. 
	It is worth noting that we are mostly interested in  $ T_\C (y^0, y^1) $,
	whereas $\overline T_\C (y^0, y^1) $ and $\widehat T_\C (y^0, y^1) $
	play an auxiliary role (although they might be of interest in their own right).
	
		\begin{comment}
	Our aim is to show that $\overline T_\C (y^0, y^1) = 	T_\C (y^0, y^1)$,
	which is a consequence of the following two lemmas.
	We use $\widehat T_\C (y^0, y^1)$ for the intermediate step,
	showing first that $\overline T_\C (y^0, y^1)= \widehat T_\C(y^0, y^1)$, and then that $\widehat T_\C(y^0, y^1) = T_\C(y^0, y^1)$.

	We start with the following technical lemma 
	establishing the existence of 
	an `inverse' map to $A_H$.
	
	\begin{lem}\label{chafe}
		There exists  a continuous map $f:A_H (H^\perp) \to H^\perp$ such that $A_H f(h) = h$
		for all $ h \in A_H (H^\perp)$
	\end{lem}
	\textbf{Proof}.
	Note that $A_H (H^\perp) $ is a linear subspace of  $ H$.
	The map $A_H: H^\perp \to A_H (H^\perp)$ need not be an injection,
	however the exists a linear subspace $K \subset H^\perp$
	such that $A_H: K \to A_H (H^\perp)$ is an injection and thus bijection.
	Such a $K$ can be chosen as a maximal collection of vectors $\{h_1 ^\perp, ... ,h_{\dim(A_H (H^\perp))} ^\perp \}
	\subset H^\perp$ such that $\{A_H h_1 ^\perp, ... , A_Hh_{\dim(A_H (H^\perp))} ^\perp \}$ 
	is a collection of linearly independent vectors.
	Letting $f$ be the inverse of $A_H: K \to A_H (H^\perp)$ finishes the proof. 
	\qed
	\end{comment}

	\begin{lem}\label{lem hat T_C = bar T_C}
		 Let $y^0, y^1 \in C ^{\mathrm{o}}$. It holds that
		\begin{equation}
		\begin{split}
		\widehat T_\C (y^0, y^1) = \overline T_\C (y^0, y^1).
		\end{split}
		\end{equation}
	\end{lem}
	\textbf{Proof}. Since the infimum  taken over a smaller set is larger,  $\widehat T_\C (y^0, y^1) \geq \overline T_\C (y^0, y^1)$.
	Let $\varepsilon >0$ be a small number. There exists $T \leq \overline T_\C (y^0, y^1) + \varepsilon$
	and $z(t) \in C([0,T], H)$ such that  
	\[
	z(0) = P_{H }y^{0}, \ \ 
	z(T) = P_{H }y^{1}, 
	\]
	and
   \eqref{charade} holds.
	%It follows from \eqref{charade} that $\dot z(t) - 	A_{H }z(t) \in A _H (H ^\perp)$.
	Set 
	\begin{equation}
	g ^ \perp (t) := f(\dot z(t), z(t)),
	\end{equation}
	where $f$ is the map from Condition \ref{con measurable selection}. Note that 
	 $\dot z(t): [0,T] \to H $ is Lebesgue measurable, and hence
	$g: [0,T] \to H ^\perp$ is Lebesgue measurable as well. Also,
	\[
	\dot z(t) = F_{H }(z(t) + f(\dot z (t), z(t)) ) =  F_{H }(z(t) + g ^ \perp (t) )  ,
	\]
	and hence $\widehat T_\C (y^0, y^1) \leq T \leq \overline T_\C (y^0, y^1) + \varepsilon$.
	Since $\varepsilon >0$ is arbitrary, the proof is complete.
	\qed

	\begin{lem}\label{T_C leq bar T_C}
	Let $y^0, y^1 \in C ^{\mathrm{o}}$. It holds that
	\begin{equation}
	\begin{split}
	\overline T_\C (y^0, y^1) \leq  T_\C (y^0, y^1).
	\end{split}
	\end{equation}
	\end{lem}
  \textbf{Proof}.  	The   statement is a consequence of 
  the fact that  the infimum 
  taken over a larger set is smaller.
  \qed
  
	\begin{rmk}
		Lemma \ref{lem hat T_C = bar T_C} and Lemma \ref{T_C leq bar T_C}
		do not require Condition \ref{con on C}:
		their conclusions hold for arbitrary measurable $\C$.
	\end{rmk}

	For $h, h_0, h_1 \in H$, $h_0 \ne h_1$, we define 
$$s(h, h_0, h_1) = \frac{\sup\langle F_H ((h + H^\perp)\cap \C) ,h_1 - h_0 \rangle  }{|h_1 - h_0|},$$
where 
\begin{equation}
\langle F_H ((h + H^\perp)\cap \C) ,h_1 - h_0 \rangle 
%\\
= \{ \langle x ,h_1 - h_0 \rangle \mid x \in F_H ((h + H^\perp)\cap \C)  \}.
\end{equation}

	The next theorem gives  a way to compute $T_\C (y^0, y^1) $
	in the case $m = n - 1$, i.e. when the range of $B$
	has co-dimension one.

	\begin{thm}\label{thm dim H = 1}
		Let $\dim H = 1$ and $h_0 \ne h_1$, where $h_i = P_H y^i$, $i = 1,2$. 
				\begin{itemize}
			\item[]$(i)$ 
			if  $s(h, h_0, h_1) > 0 $
			for all $h \in [h_0, h_1]$,
			then
				 $T_\C (y^0, y^1) = \overline T_\C (y^0, y^1) < \infty$,
			and 
			\begin{equation} \label{time formula}
			T_\C (y^0, y^1) = \int\limits _{h \in [h_0, h_1]} \frac{dh}{s(h, h_0, h_1)}.
			\end{equation}
			In particular, the integral in \eqref{time formula} is finite.
			\item[]$(ii)$ if $s(h_2, h_0, h_1) \leq 0 $ 
			for some $h_2 \in [h_0, h_1]$, 
			then $T_\C (y^0, y^1) = \overline T_\C (y^0, y^1) = \infty$.
		\end{itemize}
	\end{thm}
	\textbf{Proof}.
	We start with $(i)$.
	We begin with the auxiliary claim 
				\begin{equation} \label{time formula overline}
	\overline T_\C (y^0, y^1) = \int\limits _{h \in [h_0, h_1]} \frac{dh}{s(h, h_0, h_1)}.
	\end{equation}
	
	 Let  $s(h, h_0, h_1) > 0 $ for all  $h \in [h_0, h_1]$.
	 By definition of $s(h, h_0, h_1) $ and since $F$
	 is continuous, we have in fact $ \inf\limits_{h \in [h_0, h_1]} s(h, h_0, h_1) >0 $.
	 Hence the integral in \eqref{time formula} is well defined and finite.
	 For a small positive $\delta < \frac 12 \inf\limits_{h \in [h_0, h_1]} s(h, h_0, h_1) $, 
	the  system 
		\begin{equation}\label{pungent}
	\begin{cases}
	z(0) = P_{H}y^{0}, \\
	z(T) = P_{H}y^{1}, \\
	\dot z(t) \in F_{H }\left((z(t) + H^\perp) \cap \C \right), \\
	T >0
	\end{cases}
	\end{equation}
	 has a solution satisfying 
	$\frac{\langle \dot z(t), h_1 - h_0 \rangle}{|h_1 - h_0|} > s(z(t), h_0, h_1) - \delta$  with 
	\begin{equation}
	T^{(\delta )} = \int\limits _{h \in [h_0, h_1]} \frac{dh}{s(h, h_0, h_1) - \delta}.
	\end{equation}
	Taking  the limit $\delta \downarrow 0$, we get by the dominated convergence theorem 
	\begin{equation}
	 \overline T_\C (y^0, y^1) \leq \liminf\limits _{\delta \downarrow 0 }T^{(\delta )} = \int\limits _{h \in [h_0, h_1]} \frac{dh}{s(h, h_0, h_1)}.
	\end{equation}

	Using  Lemma \ref{speedBound} in the Appendix,
		we now prove 
	the reverse inequality $\overline T_\C (y^0, y^1) \geq  \int\limits _{h \in [h_0, h_1]} \frac{dh}{s(h, h_0, h_1)}$.
	Note that
	\begin{equation*}
	\begin{gathered}
	\frac{\langle \dot z(t), h_1 - h_0  \rangle}{|h_1 - h_0|}  \leq 
	\sup  \left\langle F_{H }\left((z(t)   + H ^\perp)\right), h_1 - h_0  \right\rangle  
	\leq   { s}(z(t), h_0, h_1)
	\end{gathered}
	\end{equation*}
	and 
	$\frac{\langle \dot z(T) - h_0, h_1 - h_0  \rangle}{|h_1 - h_0|} = |h_1 - h_0|$.
    Applying Lemma \ref{speedBound} with 
	$M  = |h_1 - h_0|$,
	$f(t) = \frac{\langle z(t) - h_0, h_1 - h_0  \rangle}{|h_1 - h_0|} $,
	and 
		 $$
	g(v) = s\left(h_0 + \frac{v}{|h_1 - h_0|} (h_1 - h_0) , h_0, h_1\right),  \ \ \ v \in [0, |h_1 - h_0|],
	$$ 
		we get
	\begin{equation}
	\overline T_\C (y^0, y^1) \geq 
	\int\limits _{0} ^{|h_1 - h_0|} \frac{dv}{ s(h_0 + \frac{v}{|h_1 - h_0|} (h_1 - h_0) , h_0, h_1)} =
	\int\limits _{h \in [h_0, h_1]} \frac{dh}{s(h, h_0, h_1)}.
	\end{equation}

	Thus, \eqref{time formula overline} is proved.
	Next we proceed with the proof of $T_\C (y^0, y^1) = \overline T_\C (y^0, y^1) $.	
	Recall that by Lemma \ref{lem hat T_C = bar T_C}, $\widehat{T} (y^0, y^1) = \overline{T} (y^0, y^1)$.
	Take $\varepsilon > 0$. There exist $T < \widehat{T} (y^0, y^1) + \varepsilon$,
	$z \in C([0,T], H )$, and
	$h^\perp \in L^0([0,T], H ^\perp )$
	such that \eqref{system hat} holds,
	$P_H y^0 = z(0)$, and $P_H y^1 = z(T)$.
	Take now a small $\delta >0$.
	Since $\C$ is bounded and convex,
	it is possible to choose 
	$h^c \in C^1([0,T], H ^\perp )$ (the space of continuously differentiable
	functions) in such a way that 
	$|h^\perp  - h^c|_{L ^1 ([0,T ], H ^\perp) } <\delta$
	and for $t \geq 0$,
	$z^c(t) +h^c(t) \in \C $,
	where $z ^c (t)$ is the	the solution to 
	\begin{equation}\label{winnow}
	\dot z ^c (t) = F_{H }(z^c(t) + h^c(t)  ), \ \ \ z ^c(0) = P_{H}y^{0}.
	\end{equation}

 	Let $y^c (t) = z ^c(t) +  h^c(t)$. 
 	Note that $y^c$ is a solution to \eqref{the model}. 	
 	Subtracting \eqref{winnow} from \eqref{system hat} we get in the integral form
 	
 	\begin{equation*}
 	z(t)  -  z ^c (t)  
 	= \int\limits _0 ^{ { t}}  F_H(z(t) - z^c (t)) dt +  \int\limits _0 ^{ { t}} F_H(h ^\perp(t) - h^c (t)) dt,
 	\end{equation*}
 	hence 
 	\begin{equation}\label{icky}
 	|z(t) - z^c (t)| \leq L_F \int\limits _0 ^t |z(t) - z^c (t)| dt + { L_F} T \delta,
 	\end{equation}
 	where $L_F$ is the Lipschitz constant for  $F$.
 	By Gr\"onwall's inequality from \eqref{icky} we obtain
 	\begin{equation}
 	|z(t) - z^c (t)| \leq L_F T e^{L_F T } \delta,  \ \ \ t \in [0,T]. 
 	\end{equation}
 	
 	In particular, 
 	\begin{equation} \label{obtuse}
 	{|h_1 - z^c (T)| }  = |z(T) - z^c (T)| \leq L_F T e^{L_F T } \delta
 	\end{equation} 
 	Recall that we took $T < \widehat{T} (y^0, y^1) + \varepsilon$.
 		Since $\C$ is compact and convex, 
 		$\inf\limits _{h \in [h_0, h_1]} {s(h, h_0, h_1)} > 0$.
 		Let $s _m >0 $ be such that 
 	\begin{equation}
 	s(h, h_0, h_1) > s _m, \ \ \ h \in [h_0, h_1].
 	\end{equation}

 	Without loss of generality we can  assume that 
 	\begin{align}
 	\langle z^c (T)  - h_0 , h_1 - h_0 \rangle & < \langle  h_1 - h_0 , h_1 - h_0 \rangle, \label{sighted1}
 	\\
 	\ \ \ \langle F(z^c (T) +  h^c(T)) , h_1 - h_0 \rangle & = s _m  > 0 \label{sighted2}
 	 	\end{align} 
 	and that $\delta >0 $ is so small  that
	there exists $r > |h_1 - z^c (T)| $ such that
	for all $y $ with $|y -z^c (T) -  h^c(T) | < r$,
	\begin{equation}\label{poppy}
	\langle F(y) , h_1 - h_0 \rangle > \frac 12 s _m
	\end{equation}
	and the ball
    	\begin{equation}\label{supple}
    \{y \in \R ^n \mid |y-z^c (T) -  h^c(T)|  \leq r \} \subset \C ^{\mathrm{o}}.
    \end{equation}  
	Indeed, if the first inequality in \eqref{sighted1} does not hold, 
	then we can go back in time. More precisely,
	we can
	replace $z^c(T)$ with $z^c(T - \Delta)$
	for a small $\Delta >0$
	so that both \eqref{obtuse} and \eqref{sighted1} hold; then,
	since $s(h, h_0, h_1)>0$ for $h \in [h_0, h_1]$,
	for some  $h^{\perp,2} $ with $z^c(T) + h^{\perp,2} \in \C$ we have
	$\langle F(z^c (T) +  h^{\perp,2}) , h_1 - h_0 \rangle { > 0}$,
	and therefore we can make 
	 the  inequality 
	in \eqref{sighted2} hold true as well 
	by modifying if necessary $h^c$ near 
	$T - \Delta$ to ensure $h^c(T - \Delta ) = h^{\perp,2} $.
Finally, \eqref{poppy} and \eqref{supple} are possible because $\C$ is  compact
and convex
and $F$ is uniformly continuous on $\C$.

	\begin{wrapfigure}{r}{0.5\textwidth}
		\vspace{-1.0cm}
		\begin{center}
			\includegraphics[width=0.48\textwidth,height=0.48\textwidth]{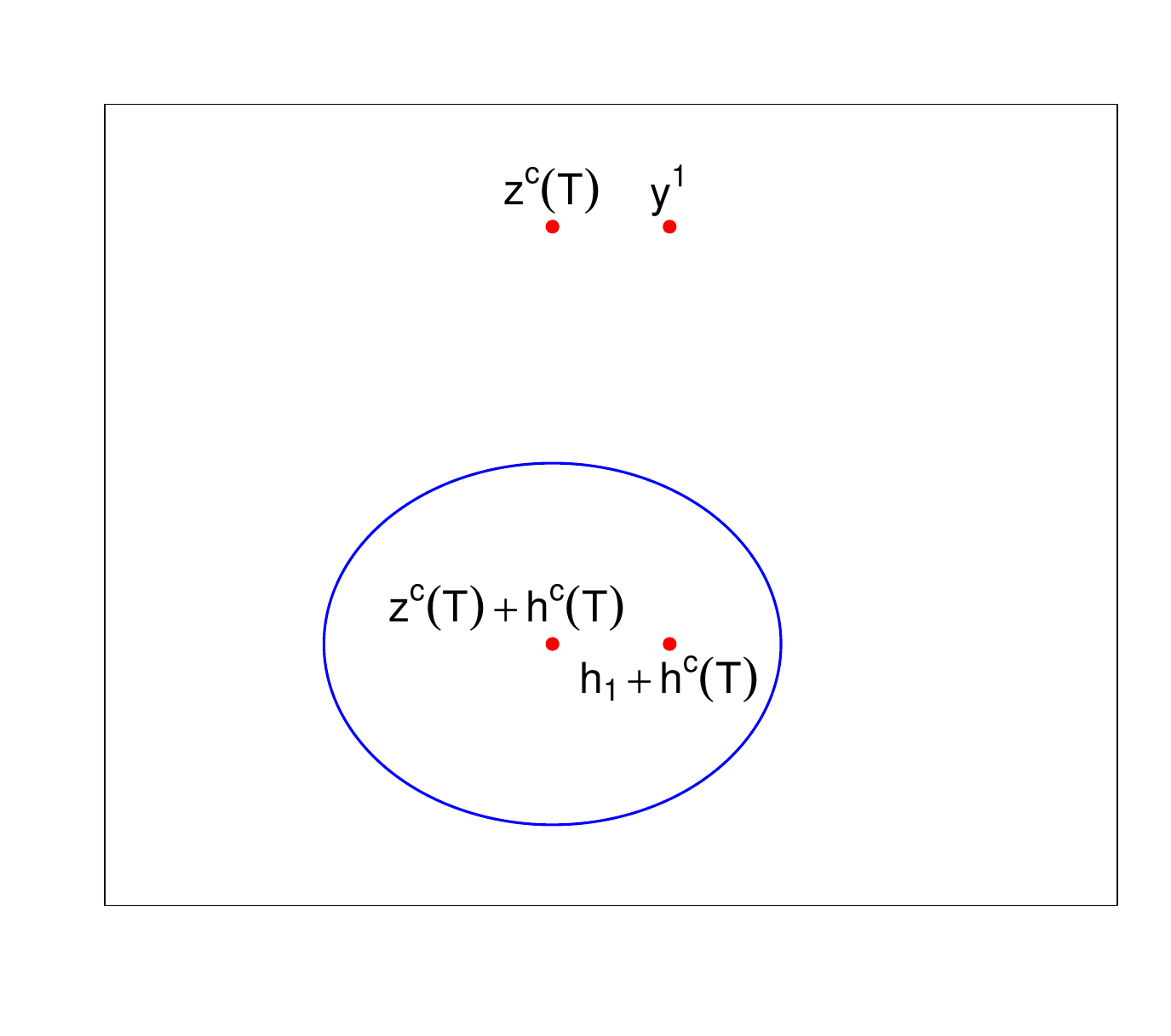}
		\end{center}
		\vspace{-1cm}
		\caption{The plane $\mathcal{L}$. The radius of the blue circle is $r$. }	\label{f-2}
		\vspace{-0pt}
		
	\end{wrapfigure}
	
	It is possible to reach $y^1$ starting from 
	$z^c (T) +  h^c(T)$ in a short time interval $[T, T_1]$.
	Indeed,
	denote by $\mathcal{L}$
	the 
	two-dimensional
	plane spanning points $z^c (T) +  h^c(T) $,
	$h_1 +  h^c(T)$, and $y^1$
	(the case
		 $h_1 +  h^c(T) = y^1 $ 
		   is simpler and discussed below).
		   The plane $\mathcal{L}$
		   is depicted on Figure \ref{f-2}.
		    Note that
	 $z^c (T) \in \mathcal{L}$
	since $z^c (T) = (z^c (T) +  h^c(T)) + y^1   - (h_1 +  h^c(T)) $,
		and hence also
			$H \subset \mathcal{L}$.
		Let us only consider controls $u$
		ensuring that $y(t)$ stays in $\mathcal{L}$,
		that is, $P_{\mathcal{L} ^\perp } B u(t) = - P_{\mathcal{L} ^\perp } F(y(t)) $, $t \geq T$. 
		Denote by $\mathcal{K} $ a one-dimensional subspace 
		of $\mathcal{L} $ orthogonal to $H$. 
	Starting from $t = T$, we take $P_{\mathcal{K}  } B u(t) = - P_{\mathcal{K} } F(y(t)) $
	at the beginning, ensuring that $y(t)$
	is moving on the interval $[ z^c (T)+ h^c(T) , h_1+ h^c(T)]$
	toward $h_1+ h^c(T)$ at a speed at least $s_m/2$.
    
    At a time $\tau$
    when $y(t)$ is near  $h_1+ h^c(T)$, we 
    stop requiring  $P_{\mathcal{K}  } B u(t) = - P_{\mathcal{K} } F(y(t)) $
    and instead take $P_{\mathcal{K}  } B u(t) = M (y^1 - h_1 - h^c(T)) $
    for a  large number $M$.
    By an intermediate value theorem, for some $\tau>T$,
     $y(t)$ is going to pass through $y^1$
     at a certain time $T_1$. 
     By taking $\delta$ small and $M$ large, we can ensure that $T_1 - T$
     is small.
     Note that here it is important that $y^1$ is in the interior $ \C ^{\mathrm{o}}$, because
     otherwise even for large $M >0$,
     a trajectory of $y(t)$ hitting $y^1$ 
      may cross the boundary of $\C$, violating the constraint condition.

     In the case  $h_1 +  h^c(T) = y^1 $ we 
     just
     require  $P_{H ^\perp } B u(t) = - P_{H ^\perp } F(y(t)) $
     for $t \geq T$, ensuring that $y(t)$
     stays on the interval $[z^c (T) +  h^c(T), h_1 +  h^c(T) ] = [z^c (T) +  h^c(T), y^1 ]$
     and hits $y^1$.

     We note here that in particular if $ \langle P _H F (y^1), h_1 - h_0 \rangle < 0  $,
     it is important that $y ^1 \in \C ^\mathrm{o}$, because the trajectory
     of $y(t)$ has to reach $y^1$ from `behind', that is, from within the 
     half-space $ \left\{ u \in \R ^n: \langle u - h_1 , h_1 - h_0 \rangle > 0  \right\} $.

	Therefore, $T_\C (y^0, y^1) = \widehat T_\C (y^0, y^1) =\overline T_\C (y^0, y^1) $,
	and $(i)$ follows from \eqref{time formula overline}.

	Having proved $(i)$, we now turn to $(ii)$. Let $h_2 \in [h_0, h_1]$
	be such that $s(h_2, h_0, h_1) \leq 0 $.
	Since the boundary $\partial \C$ is assumed to be differentiable
	and $\C$ is convex, the Borel set valued map
	\begin{equation}
	 [h_0, h_1] \ni h \mapsto \C _h \in \mathscr{B}(H^\perp)
	\end{equation}
	defined by $ \C _h = \{ y \in \C\mid P_H y  = h \}$
	is Lipschitz continuous
	in Hausdorff distance (or Hausdorff metric, see e.g. \cite{Hausdorff_metric})
	 with some constant $L_1$.
	Since $F$ is also Lipschitz continuous,
	it  holds that 
	\begin{equation}\label{wean}
	s(h, h_0, h_1) \leq (L_F + L_1)|h - h_2|, \ \ \ h \in  [h_0, h_2].
	\end{equation}
    It follows from \eqref{wean} that any solution to \eqref{charade}
    starting from $P_H y^0$
    never reaches the location $h_2$:
    that is,
    for all $t \geq 0 $,
    $$\langle z(t) - h_0, h_1 - h_0 \rangle <  \langle h_2 - h_0, h_1 - h_0 \rangle. $$
	Consequently, \eqref{pungent} has no solutions. Hence
	$ \overline T_\C (y^0, y^1) = \infty$,
	and by Lemma \ref{T_C leq bar T_C}, 
		$T_\C (y^0, y^1)  = \infty$.
		\qed

   The purpose of following example is to demonstrate that
   	the assumption $\dim H = 1$
   	in
   	 Theorem \ref{thm dim H = 1} 
   is necessary: without it the conclusion of the theorem need not be true.
		
	\begin{example} \label{example bar T ne T}

	Here we provide an example 
	where all conditions of 
	Theorem \ref{thm dim H = 1},
	$(i)$ are satisfied except $\dim H = 1$,
	 but the conclusions of 
	 $(i)$ are false, in particular,
	 $$
	 \overline T_{\C } (y^0, y^1) = \widehat T_{\C } (y^0, y^1) \ne  T_{\C } (y^0, y^1) .
	 $$ 
	Let 
	$n = 3$, $m= 1$, $\C = [-2,2] ^3$, $B : \R \to \R ^3$,
	$B u = (0,0,u)^\top$.
	Thus, $H = \{(\alpha,\beta,0)^\top\mid \alpha, \beta \in \R  \}$,
	and $H ^\perp = \{(0,0,\gamma)^\top\mid  \gamma \in \R  \}$.
	Note that $\dim H = 2$.
	We impose the following conditions on $F$.
	\begin{enumerate}
		\item For $x \in \C$, $x \notin
		 \Big\{(\alpha, 0, 1 )\mid \alpha \in [-1,1] \Big\} \cup \Big\{(\alpha, 0, 0 )\mid \alpha\in [-1,1]  \Big\}  $,
		we have
		$\langle F(x), (0,1,0)^\top \rangle > 0$. Here and elsewhere, $\langle \cdot , \cdot \rangle$ is the scalar product in $\R ^n$.
		\item $F(\alpha, 0, 1)^\top = F(\beta, 0, 0)^\top =  (1,0,0)^\top$ for $\alpha \in [-1,0]$,
		$\beta \in [0,1]$.
		\item    $F(\alpha, 0, 1)^\top = F(\beta, 0, 0)^\top =  (-1,0,0)^\top$ for $\alpha \in [\frac 12,1]$,
		$\beta \in [-1, -\frac 12]$.
	\end{enumerate}
	Of course,  many vector fields exist satisfying those conditions.	
	Take now $y^0 = (-1,0,1)^\top$ and $y^1 = (1,0,0)^\top$. It is not difficult to see that 
	$T_\C (y^0, y^1) = \infty$, since once the trajectory $y(t)$
	of the solution to \eqref{the model} leaves the segment $\{(\alpha, 0, 1 )\mid \alpha \in [-1,1] \}$,
	the second coordinate of $y(t)$ becomes positive and stays positive forever:
	$ \langle y(t) , (0,1,0)^\top \rangle >0$.
	On the other hand, $\widehat T_\C (y^0, y^1) \leq 2$.
	Indeed, 
	$z(t) = (-1 +t, 0, 0)^\top \in H$ 
	and 
	$$h^\perp (t) = \begin{cases}
	(0,0, 1)^\top,  \ \ t < 1, \\
	(0,0, 0) ^\top, \ \ t \geq 1,
	\end{cases}
	$$
	give a  solution to \eqref{system hat} with $z(0) = P_H y^0$
	and
	$z(2) = P_H y^1$.
	Thus, $ \widehat T_{\C } (y^0, y^1) \leq 2 < \infty =   T_{\C } (y^0, y^1) $.
	The equality 
	$ \widehat T_{\C } (y^0, y^1) = \overline 
	T_{\C } (y^0, y^1) $ follows from Lemma \ref{lem hat T_C = bar T_C}.
	The example can of course be modified so that $T_{\C } (y^0, y^1) < \infty$
	but still $T_{\C } (y^0, y^1) > 2 \geq \widehat T_{\C } (y^0, y^1)$.
			\hfill
	$\Diamond$
		\end{example}

\begin{figure}[H]
		\vspace{-0.6cm}
		\begin{center}
		\includegraphics[scale=0.49]{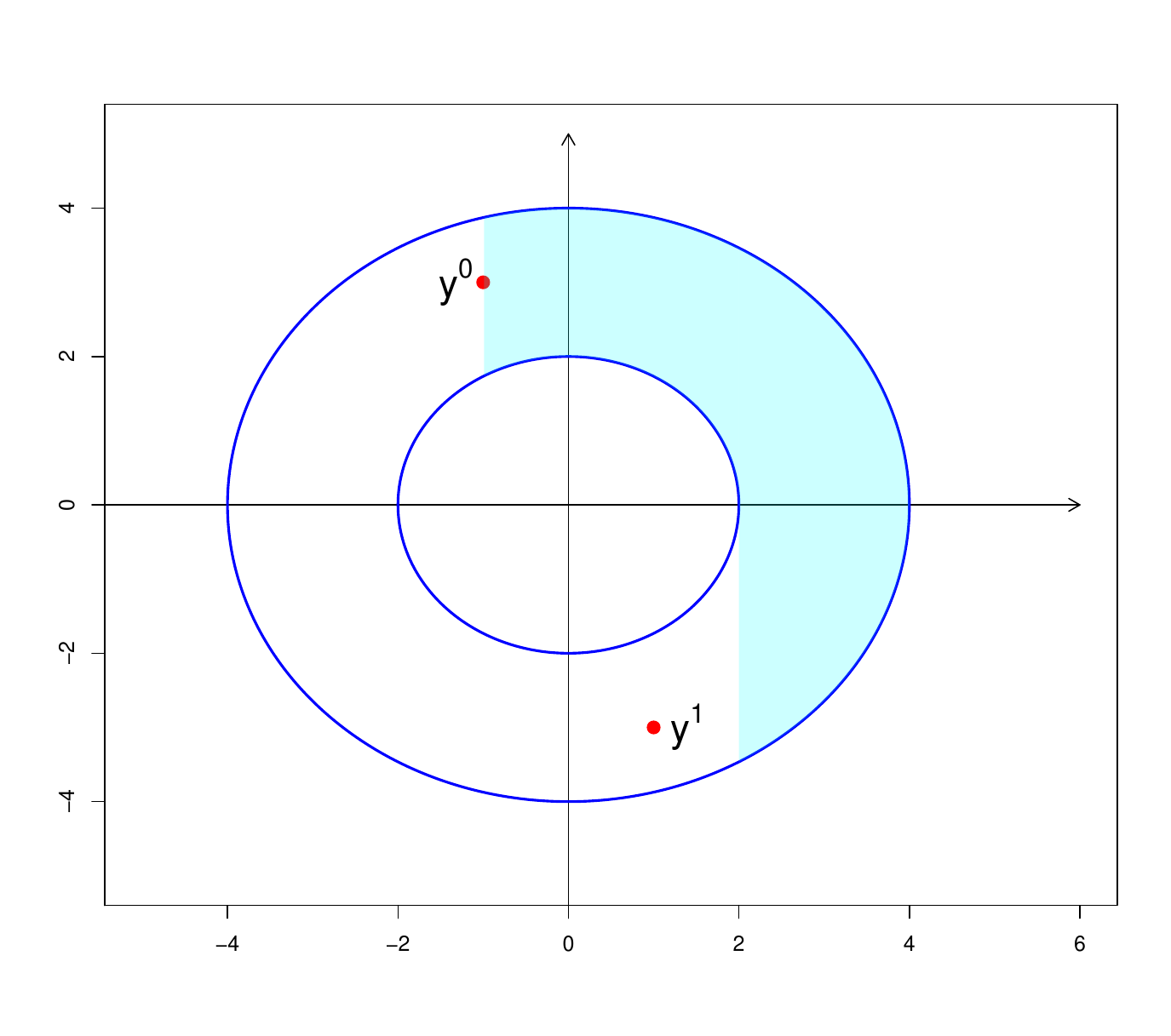}
	\end{center}
\vspace{-1.3cm}
		\caption{An illustration to Example \ref{example with annulus}. 
		The boundaries of the annulus are blue. The shaded area are 
		the points reachable under the 
		  state constraint $\C$. }
		\label{f-1}
			\vspace{-0.6cm}
	\end{figure}

	\begin{example} \label{example with annulus}
		Let us also mention an annulus as an example where all conditions 
		of Theorem \ref{thm dim H = 1} are satisfied
		except the convexity of $\C$, while the conclusions of Theorem \ref{thm dim H = 1} fail. Take $n = 2$, $F( (x_1, x_2)^{\top} ) \equiv (1
	%	{   + \frac{x_2}{10000}}
		,0)^\top$, $B = (0,1)^\top$,
		\[
		\C = \left\{ (x_1, x_2)^\top \middle| 4 \leq x_1 ^2 + x_2 ^2 \leq 16  \right\},
		\]
			$y^0 = (-1,3)^\top$, $y^1 = (1,-3)^\top$ (see Figure \ref{f-1}). Then indeed $T_\C(y^0, y^1) = \infty$
			since starting from $y^0$ it is not possible to reach the part of the annulus below the hole.
					\hfill
			$\Diamond$
	\end{example}

	%\vspace*{0.4cm}

	We now introduce another reachability time.
	To start off, we define a solution 
	satisfying the property that 
	a small perturbation at any point
	does not break the reachability property.
	In the definition below we use the terms
	`reachability' or `reachable' with regard to system \eqref{the model}.
	\begin{defi} \label{iniquity}

	Take $y^0, y^1 \in \C ^o$,
	and let $z(t), h^\perp(t)$ be a solution of \eqref{system hat}
	with $z(0)  =P_H y^0 $, $z(T)  =P_H y^1 $. 
	We call
	this solution \emph{pliable}
	if
	\begin{itemize}
		\item []$(i)$ 	for every $t_1$ and $ t_2$, $0\leq t_1 < t_2 \leq T$,
		there exists $\varepsilon _{t_1, t_2} >0$
		and $\tau(t_1, t_2) > 0$
		with the property that 
		for all $ \varepsilon \leq \varepsilon_{t_1, t_2}$
		there is $\delta _{\varepsilon, t_1, t_2} >0$
		satisfying the following condition:
		for every $y_{t_1} \in \mathcal{B}(z(t_1) + h^\perp(t_1) , \delta _{\varepsilon, t_1, t_2} ){ \cap \C ^{\mathrm{o}}}$
		there exists $ y_{t_2} \in  \mathcal{B}(z(t_2) + h^\perp(t_2) , \varepsilon ){ \cap \C ^{\mathrm{o}}}$
		such that 
		$y_{t_2}$ is reachable from $ y_{t_1}$
		with state constraint $\C$ within the time  $\tau(t_1, t_2)$ (that is,
		in a time not greater than $\tau(t_1, t_2)$).
		\item[]$(ii)$ For all $t \in [0,T]$ sufficiently close to $T$
		there exists small $\delta _t >0$
		such that $y^1$ is reachable from every point 
		of $\mathcal{B}(z(t) + h^\perp(t) , \delta _t )$ with the state constraint $\C$,
		and
		\begin{equation}
		 \limsup\limits _{t \uparrow T } \sup\limits_{y \in \mathcal{B}(z(t) + h^\perp(t) , \delta _t ){ \cap \C ^{\mathrm{o}}}} T_{\C}(y, y^1) = 0.
		\end{equation}
		\item[] $(iii)$ 	The following inequality holds
		\begin{equation}\label{ancillary}
		\limsup\limits _{\Delta t \to 0} \sup \limits _{t \in [0, T - \Delta t]} \frac{\tau(t, t + \Delta t)}{\Delta t} \leq 1. 
		\end{equation}

	\end{itemize}
	\end{defi}

	While the definition of a pliable solution
	may seem unwieldy,
	intuitively it means that a solution to \eqref{system hat}
	can be approximate well by solutions to \eqref{the model}
	along the entire trajectory
	without incurring significant  loss of time.

	Now we define
	\begin{equation}
	\begin{split}
	\widetilde T_\C (y^0, y^1) =  \label{widetilde T def}
	\inf\big\{T: \text{there exists } z(t) \in C([0,T], H ),
	h^\perp(t) \in L^0([0,T], H ^\perp )
	\text{ s.t. } 
	z(0) = P_{H}y^{0},
	\\
	z(T) = P_{H}y^{1}, 
	\text{ \eqref{system hat} holds and the solution is pliable}  \big\}. \notag
	\end{split}
	\end{equation}
	
	The following result gives an upper bound of the reachability time $T_\C (y^0, y^1)$.
	
	\begin{prop}\label{T leq tilde T}
		It holds that 
		\begin{equation}
		  T_\C (y^0, y^1) \leq \widetilde T_\C (y^0, y^1).
		\end{equation}
	\end{prop}
	\textbf{Proof}. Let
	$\varepsilon >0$ and let 
	 $z(t) \in C([0,T], H )$
	and
	$h^\perp(t) \in L^0([0,T], H ^\perp )$ constitute a pliable solution
	to \eqref{system hat} with $z(0) = P_{H}y^{0}$,
	$z(T) = P_{H}y^{1}$, $T \leq \widetilde T_\C (y^0, y^1) + \varepsilon$.
	Take a partition $0 = t_0 < t_1 < ... < t_k < T $
	such that
	 \begin{equation}\label{obtuse angle}
	 \sum\limits_{i = 0} ^{k-1} \tau (t_i, t_{i+1}) \leq T + \varepsilon,
	\end{equation}
	and for some $\delta _k >0$
		 \begin{equation}\label{incarnate}
	 \sup\limits_{y \in \mathcal{B}(z(t_k) + h^\perp(t_k) , \delta _k )} T_{\C}(y, y^1)  \leq  \varepsilon.
	\end{equation}
	We note that \eqref{obtuse angle} 
	and 
	\eqref{incarnate}
	are possible  by items $(ii)$
	and $(iii)$ 
	of Definition \ref{iniquity}, respectively.
	Next we define the sequence $r_k , r _{k-1}, ..., r_1$
	consecutively as follows:
	$r_k = \min(  \delta _k, \varepsilon_{t_{k-1}, t_{k}})$,
	$r_{k-1} = \min ( \delta _{r_k, t_{k-1}, t_k},  \varepsilon _{t_{k-2}, t_{k-1}}  )$,
	$r_{k-2} = \min ( \delta _{r_{k-1}, t_{k-2}, t_{k-1}},  \varepsilon _{t_{k-3}, t_{k-2}}  )$
	and so on, until 
	$r_{1} = \min ( \delta _{r_{2}, t_{1}, t_{2}},  \varepsilon _{t_{0}, t_{1}}  )$.
	It follows from Definition \ref{iniquity}
	that there is a solution $y(t)$ to \eqref{the model}
	satisfying 	
	$$y\left(\sum\limits_{i = 0} ^{j-1} \tau (t_i, t_{i+1}) \right) - (z(t_j) +h^\perp(t_j) ) \leq r_k,
	\ \ \ j = 1,2,..., k-2,$$
	and 
	$
	y\left(\sum\limits_{i = 0} ^{k-1} \tau (t_i, t_{i+1}) \right) \in \mathcal{B}(z(t_k) + h^\perp(t_k) , \delta _k )
	$.
	It follows from \eqref{incarnate} 
	that $y(t)$
	can be extended to reach $y^1$
	by the time $\sum\limits_{i = 0} ^{k-1} \tau (t_i, t_{i+1}) + \varepsilon $,
	that is, by the time $ T + 2 \varepsilon$ if we take \eqref{obtuse angle} into account.
	Since $\varepsilon >0$ is arbitrary, this completes the proof. \qed
	
	Revisiting Example \ref{example bar T ne T},
		 we see by continuity of $F$ 
	that the solution to \eqref{system hat}
	given there satisfies $(i)$ and $(iii)$
	of Definition \ref{iniquity}, but does not satisfy $(ii)$.
	Thus, if we modified $F$ near $y^1 = (1,0,0)^\top $
	in such a way that $(ii)$ 
		of Definition \ref{iniquity}
	was satisfied,
	then by Proposition \ref{T leq tilde T}
	we would have $T_\C (y^0, y^1) = 2 $.

	In the next theorem we give a lower bound
	on the reachability time in the case $m< n-1$.
	
	\begin{thm} 
		Let $m < n-1$. Denote $h_k = P_H y^k$, $k =0,1$.
		Denote also by $H_\perp$ the subspace of $H$ of co-dimension one
		such that $ H_\perp \perp h_1 - h_0$.
		Define
		$$\bar s(h, h_0, h_1)  = \frac{\max\langle F ((h + H_\perp + H^\perp)\cap \C) , h_1 - h_0\rangle  }{|h_1 - h_0|}.$$
		Then 
		\begin{equation}
		T_\C (y^0, y^1) \geq  \int\limits _{h \in [h_1, h_2]} \frac{dh}{\bar s(h, h_0, h_1)}.
		\end{equation}
		
	\end{thm}
	
	\textbf{Proof}. 
	Let $y(t)$, $t \in [0,T]$, be the a solution to \eqref{the model} with constraint $y(t) \in \C$.
	Let $x(t) $ be the orthogonal projection of $y(t)$ on the line spanned by $h_1 - h_0$.
	Then $x(0) = h_0$, $x(T) = h_1$, and
	\begin{equation*}
	\dot x(t) \in F_{H }\left((x(t) + H_\perp + H ^\perp) \cap \C \right).
	\end{equation*}
	Hence 
	\begin{equation*}
	\begin{gathered}
	\frac{\langle \dot x(t), h_1 - h_0  \rangle}{|h_1 - h_0|} 
	 \leq \sup \left\{ \left\langle F_{H }\left((x(t)  + H_\perp + H ^\perp) \cap \C \right), h_1 - h_0  \right\rangle  \right\}
	\leq \bar s(x(t), h_0, h_1).
	\end{gathered}
	\end{equation*}
	
	Since
	\[
	\int\limits_0 ^T \frac{\langle \dot x(t), h_1 - h_0  \rangle}{|h_1 - h_0|} dt = |h_1 - h_0|,
	\]
	the statement of the proposition follows from Lemma \ref{speedBound}.
	\qed

	\begin{rmk}
		
		In the case  $\overline T_\C (y^0, y^1) < \infty $,
		the signal $u(t)$ resulting in  a  time close to the infimum 
		can be computed
		as follows: let $h\in [h_0, h_1]$ be such that
		$P _H y(t) = h = h(t) $, $P _{H^\perp} y(t) = h ^\perp = h ^\perp(t) $ , $h \ne h_0, h_1$. Then $ u(t)$
		is informally determined by 
		$h, h ^\perp$ from the system
		\begin{equation} \label{optimal u}
		\begin{cases}
		\sup\{ \langle F_H   \big((h +  H^\perp)  \cap \C\big) ,h_1 - h_0 \rangle  \} 
		= \langle F_H (h + h^ \perp ) , h_1 - h_0 \rangle  
		\\
		h + h^ \perp  \in \C ,
		\\
		\dot h + \dot h ^ \perp   = F (h  + h^ \perp) + Bu.
		\end{cases}
		\end{equation}
		Note that the supremum in the first equation in \eqref{optimal u}
		is not always achieved (for example, if $\C$ is open, the maximum need not be achieved).
		In this case, we may either consider the closure of $\C$,
		or replace $\sup\{ \langle F_H   \big((h +  H^\perp)  \cap \C\big) ,h_1 - h_0 \rangle  \}$
		with some $\sup\{ \langle F_H   \big((h +  H^\perp)  \cap \C\big) ,h_1 - h_0 \rangle  \} - \delta$
		for a small $\delta >0$.
		
		After finding $u(t)$, we just set $ y(t) = h(t) + B u(t)$.
	We note that, typically, the infimum in \eqref{def T_C}
	can not be achieved with an $L^\infty$ control, see \cite{LTZ18}
	for the linear case.
	Most of the time system \eqref{optimal u}
	would have a solution only if
	we allow the  impulse control,
	i.e., we would allow 
	$u$ to take value in some space of distribution.
	This is the approach taken in \cite{RadonControls}
	for a linear system with a control constraint.
	In this way we could handle an	
	 instantaneous movement along a  direction from $H^\perp$.
	 To stay within $L^\infty$ controls,
	  we may need to approximate the solution 
	  to \eqref{optimal u} with $L^\infty$ controls.
	  The approximation should be possible in most cases,
	  however care needs to be taken to do the approximation properly,
	  and it is impossible to do in certain situations
	  as demonstrated by Example \ref{example bar T ne T}.
	  The approximation is  discussed in \cite{RadonControls}.
	
	 \begin{comment}
	 In our settings the infimum should be  achievable  with
	an  impulse control, since $\C$ is closed.
	To handle the latter we would need the control to take 
	value in some space of distributions.
	Using an impulse control can result in the
	`ideal'
	trajectory  satisfying \eqref{optimal u}.
	If $A_H : \R ^n \to H$ is linear, $h + h^ \perp$ from
	\eqref{optimal u} lies on the boundary of $\C$
	if it is unique.
	The `ideal' trajectory therefore is going to spend most of the time being 
	on the boundary of $\C$, while using impulse controls to move between
	different parts of the boundary if necessary.
	\end{comment}
		\end{rmk}

	\section{Examples} \label{sec examples}
	
	In this section we give three examples for which we compute the controllability time.
	In the next example we apply Theorem \ref{thm dim H = 1} to a three-dimensional control system.

	\begin{example} Let $n  = 3$, $m = 2$, $F(y) = A y$,
	
	\[
	\C = \left\{ (x_1,x_2, x_3) \mid x_3^2  + \frac{x_2^2}{4}  + \frac{x_1^2}{16}\leq 1 \right\}.
	\] 
	\[
	A = (a_{ij})_{i, j = 1,2,3} =   \begin{pmatrix}
	-0.1 &  7.0  &  0  \\
	0.0 &  0.8  &  7   \\
	0.0 & -10.0 &  -4 
	\end{pmatrix} ,
	\ \ \ 
	B =  \begin{pmatrix}
	0 & 0 \\  1 &  2 \\ 0 &  1
	\end{pmatrix} ,
	\]
	and the initial and target state
	$y^0 = (-1,0,0) ^T$, $y ^1 = (1,0,0) ^T$.
	In this example 
	$H =  \{( \kappa,0,0)^T , \kappa \in \R \}$, 
	and for $h \in [h_0, h_1]$, 
	$ h = ( \kappa,0,0)^T$, $\kappa \in [-1,1]$,
\end{example}
	\[
	s(h, h_0, h_1) 
	= \max \left\{ a_{11} \kappa + a_{12} \beta + a _{13} \gamma :  \frac{\kappa^2}{16} + \frac{\beta ^2}{4} + \gamma^2 \leq 1  \right\}
	\]
	\[
	= \max \left\{ -0.1 \kappa + 7 \beta  :   \frac{\beta ^2}{4}  \leq 1 -  \frac{\kappa^2}{16}  \right\}
	= -0.1 \kappa + 14 \sqrt{1 -  \frac{\kappa^2}{16}}.
	\]
	Hence by Theorem \ref{thm dim H = 1}
	\begin{equation}
	\begin{gathered}
	T_\C (y^0, y^1) = \int\limits _{h \in [h_0, h_1]} \frac{dh}{s(h, h_0, h_1)} 
	= 
	\int\limits _{\kappa \in [-1, 1]} \frac{d \kappa }{-0.1 \kappa + 14 \sqrt{1 -  \frac{\kappa^2}{16}}}.
	\end{gathered}
	\end{equation}
	\hfill
	$\Diamond$

	In the next example some of the points are not reachable.
	
	\begin{example} \label{mascot}
 Let $n = 2$, $m=1$,
	$\C$ be the square 
	\[
	\C = \{ (x_1,x_2) \mid  |x_1| \leq 1, |x_2| \leq 1 \},
	\] 
	and let $F(y) = A y$,
	\begin{equation}\label{A square}
	A = \begin{pmatrix}
	-2 & 3 \\ -2 & 1
	\end{pmatrix},
	%\end{equation}
	\ \ \
	%\begin{equation}
	B = \begin{pmatrix}
	1 \\  1
	\end{pmatrix},
	\end{equation}
\end{example}
	We have here $H = (\text{ran}(B))^\perp = \{ (\kappa, - \kappa) ^\top\mid \kappa \in \R \}$, 
	and

	\begin{equation}
	P_H = \begin{pmatrix}
	0.5 & -0.5 \\ -0.5 & 0.5
	\end{pmatrix},
	%\end{equation}
	%\begin{equation}
	\  \ \
	A _ H = P_H A = \begin{pmatrix}
	0 & 2 \\ 0 & -2
	\end{pmatrix}.
	\end{equation}
	
	Take 
	$y^0 =  (
	0.7 ,  -0.5
	)^\top $, $y^1 =  (
	-0.5 ,  0.3
	)^\top $. Then
	$h_0 = ( 
	0.6 ,  -0.6
	)^\top $, $h_1 = ( 
	-0.4 ,  0.4
	)^\top $.
	For $h = (\kappa, - \kappa)^T \in [h_0, h_1]$, $\kappa \in [-0.4, 0.6]$, we compute
	\[
	s(h, h_0, h_1) = \frac{\max\{\langle A_H ((h + H^\perp)\cap \C) ,h_1 - h_0 \rangle  \}}{|h_1 - h_0|}
	= 
	\frac{\max\{\langle A_H ((h + H^\perp)\cap \C) , (-1,1)^T \rangle  \}}{\sqrt{2}}
	\]
	\[
	= \frac{1}{\sqrt{2}} \max \left\{ \left\langle   \begin{pmatrix}
	0 & 2 \\ 0 & -2
	\end{pmatrix} \begin{pmatrix}
	\kappa + \alpha  \\  -\kappa + \alpha
	\end{pmatrix}  ,\begin{pmatrix}
	-1 \\  1
	\end{pmatrix}  \right\rangle \middle| 
	\begin{pmatrix}
	\kappa + \alpha  \\  -\kappa + \alpha
	\end{pmatrix} \in \C \right\}
	\]
	\[
	= \sqrt{2}  \max \Big\{ 
	2 \kappa - 2\alpha
	\Big| 
	|\kappa + \alpha| \leq 1,   |-\kappa + \alpha| \leq 1
	\Big\}  =2 \sqrt{2}( \kappa + 1 - |\kappa|),
	\]
	where we used that the maximum is achieved for $\alpha  = -1 + |\kappa|$.

	Thus, by Theorem \ref{thm dim H = 1}
	$$
	T_\C (y^0, y^1) 
	=
	\frac{1}{2\sqrt{2}}
	\int\limits _{\kappa \in [-0.4, 0.6]} \frac{d \kappa }{1 + 2\min(0,\kappa)}
	= \frac{1}{2\sqrt{2}}(0.6 + \ln 5).
	$$

\begin{figure}[H]
	
	\begin{center}
		\includegraphics[scale=0.56]{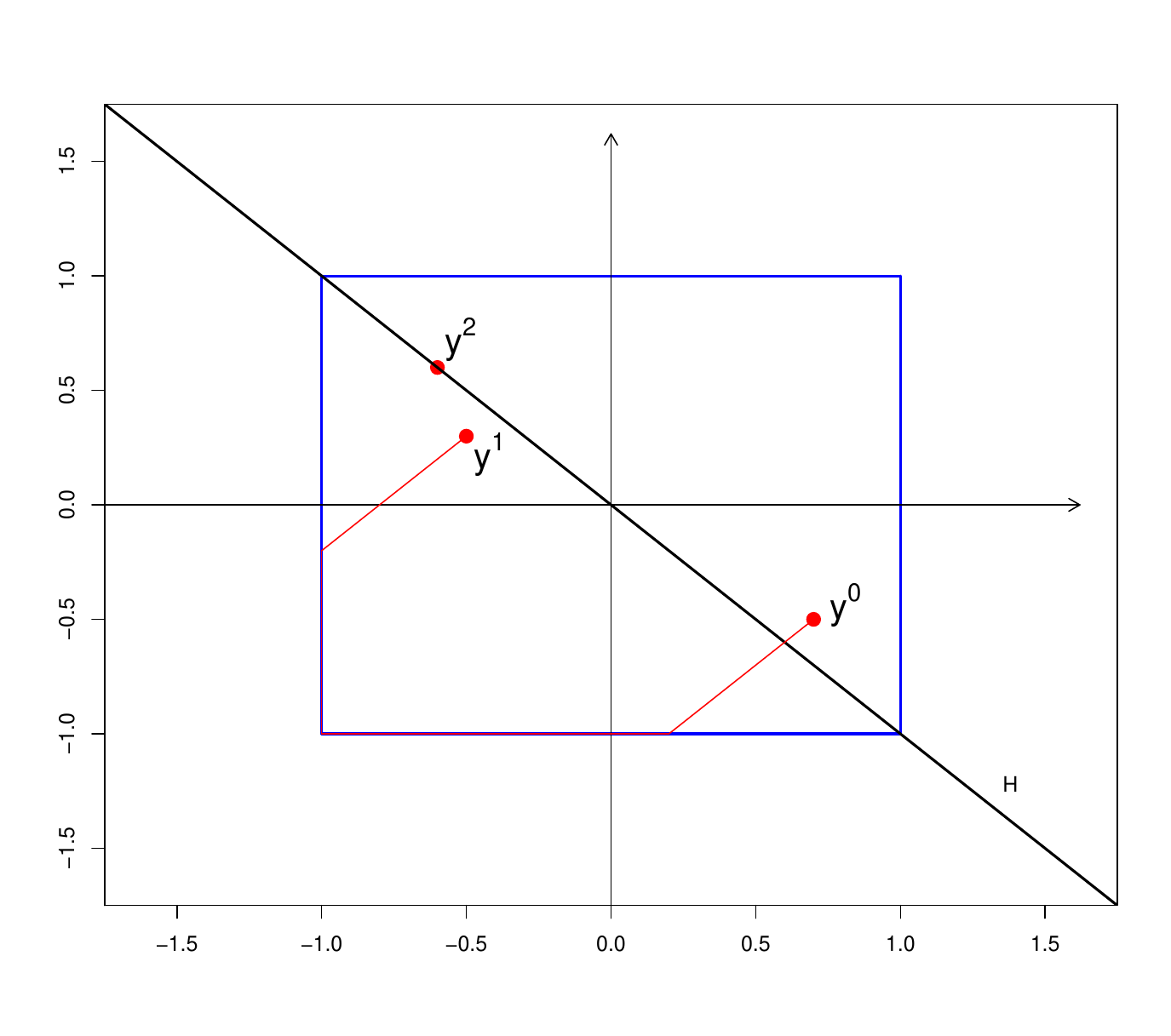}
	\end{center}
	\vspace{-1cm}
	\caption{
		The optimal  trajectory from $y^0$ to $y^1$ for \eqref{the model} in Example \ref{mascot}. The boundaries of  $\C$ are blue.
		The point $y^2$ is not reachable from $y^0$. }
	\label{f66}
	
\end{figure}

	Now, in this example not every point within $\C$ is reachable from any other.
	Take for example $y^2 = 
	(-0.6,  0.6)^T
	$, $h_2 = y^2$. Then following the same steps as above, we find that
	for
	$h = (\kappa, - \kappa)^T \in [h_0, h_2]$, $\kappa \in [-0.6, 0.6]$
	\[
	s(h, h_0, h_2) = 2 \sqrt{2}( \kappa + 1 - |\kappa|).
	\]
	Thus, $s(h, h_0, h_2) < 0 $ for example for $h = (-0.55, 0.55)^T $.
	Hence by Theorem \ref{thm dim H = 1}, 
	$
	T_\C (y^0, y^2) 
	= \infty
	$, and $y^2$ is not reachable from $y^0$.
	This example is illustrated in Figure \ref{f66}.
	\hfill
	$\Diamond$

	In the next example we deal with a non-linear system.
	
	\begin{example} \label{rustic}
		 Let $n = 2$, $m=1$, 
	\[
	\C = \{ (x_1,x_2) \mid  x_1^2 + x_2^2 \leq 5 \},
	\] 
	and let 
	$$
	F(x)
	% = \frac{1}{ \left((10 - x_1 )^2 + (4 - x_2) ^2\right)^{3/2} }
	% \left(10 - x_1,  4 - x _2\right)^\top
	 = \frac{1}{ (10 - x_1 )^2 + (4 - x_2) ^2 }
	 \frac{\left(10 - x_1,  4 - x _2\right)^\top}{|\left(10 - x_1,  4 - x _2\right)^\top|},
	 $$
	\end{example}
	$	B = \begin{pmatrix}
		0 \\  1
	\end{pmatrix}$,
	$y^0 = (
	-4 ,  -1
	)^\top$, $y^1 =( 	4 ,  -2
	)^\top$. 
	Here $F$ represents attraction of a body located at $x = (x_1, x_2)^\top$
	toward a source located at $\rm{S} = (10,4)^\top$
	with the strength of attraction being 
	reversely proportional to the square of the distance between $x$ and $\rm{S}$.
	
\begin{figure}[H]
			\vspace{-0.4cm}
		\begin{center}
		\includegraphics[scale=0.63]{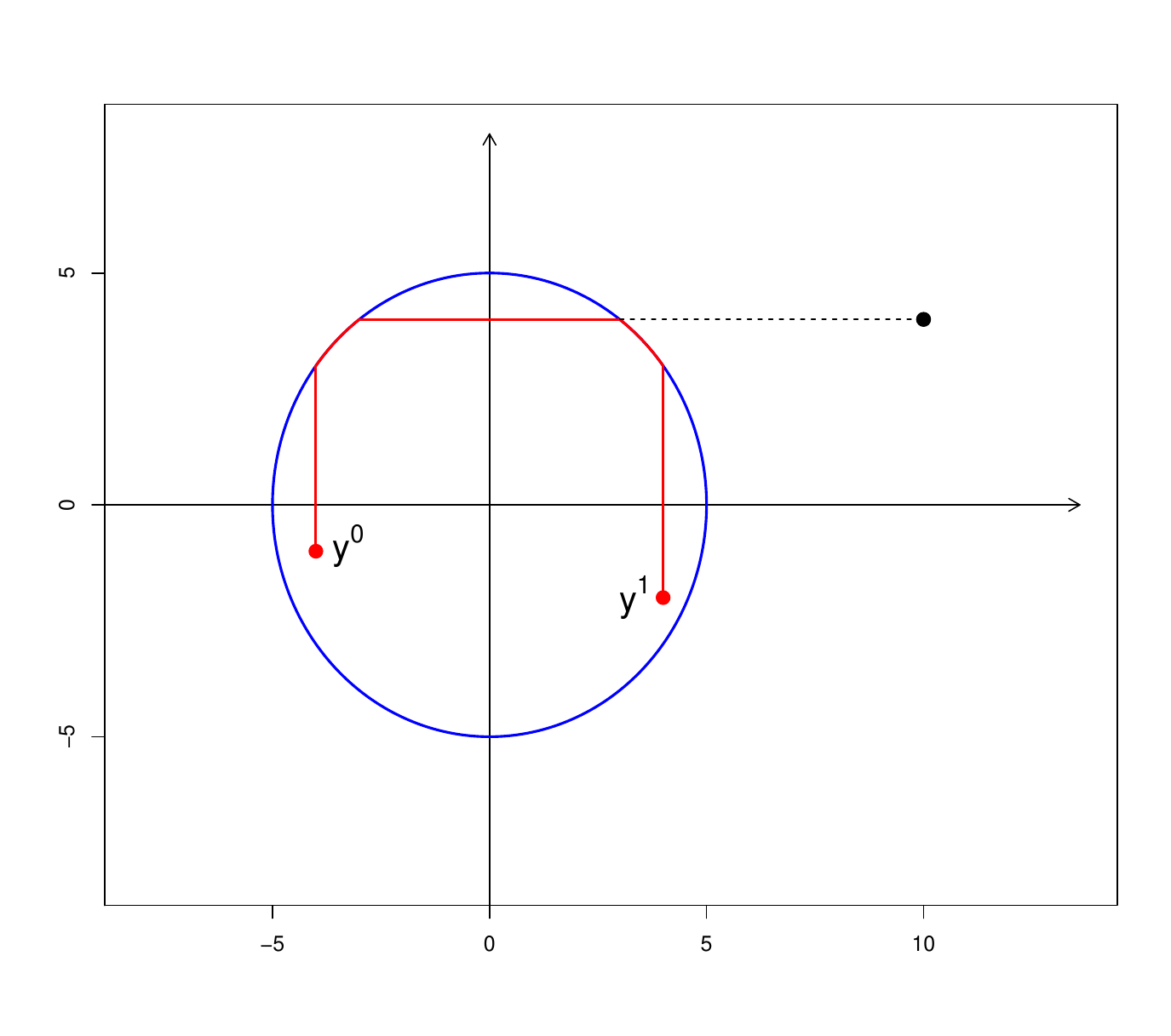}
			\end{center}
			\vspace{-1cm}
		\caption{
			The optimal  trajectory for Example \ref{rustic}. The source of attraction $\rm{S}$
		is the big black dot on the right.}
		\label{f77}
		\vspace{-1cm}
	\end{figure}

	In this case 
	$H =  \{ (\kappa, 0) ^\top \mid \kappa \in \R \}$, $h_0 = ( 
	-4 ,  0
	)^\top$,
	 $h_1 = ( 
	4 ,  0
	)^\top$. We compute
	for $h = (\kappa, 0) ^\top$
	
	\begin{equation}
	s(h, h_0, h_1) = 	\begin{cases}
	\begin{aligned}
	&\frac{ (10 - \kappa)^2 }{\big((10 - \kappa )^2 + (4 - \sqrt{25 - \kappa^2})^2\big) ^2 },& -4 \leq   \kappa  \leq 3
	\\
	&\frac{1}{(10 - \kappa )^2 },    & -3 <   \kappa  < 3,
		\\
	&\frac{ (10 - \kappa)^2 }{\big((10 - \kappa )^2 + (4 - \sqrt{25 - \kappa^2})^2\big) ^2 },& 3 \leq   \kappa  \leq 4
	\end{aligned}
	\end{cases}
	\end{equation}
	% for computation: ( (10 - x )^2 + (4 - sqrt(25 - x^2) )^2 )^2/ (10 - x )^2 
	and hence by Theorem \ref{thm dim H = 1}
	\[
	T_\C (y^0, y^1) = \int\limits _{-4 } ^{-3} \frac{\big((10 - \kappa )^2 + (4 - \sqrt{25 - \kappa^2})^2\big) ^2 }{(10 - \kappa)^2  } d \kappa
	+  \int\limits _{-3} ^3 (10 - \kappa) ^2 d \kappa 
	+  \int\limits _{3} ^{4}  \frac{\big((10 - \kappa )^2 + (4 - \sqrt{25 - \kappa^2})^2\big) ^2 }{(10 - \kappa)^2  }  d \kappa
	\]
	\[
	\approx 182.9087 +  618 + 42.9122 = 843.8209.
	\]
	This example is illustrated in Figure \ref{f77}.
		$\Diamond$

		\hfill

	\section{Conclusions and further comments} \label{sec conclusion}
	
	For a non-linear system with linear control 
	and bounded convex state constraint we 
	give results about the controllability time
	between two points.
	 The results of \cite{LTZ18}
	 are extended 
	in two directions: 
	the system
	has a non-linear drift term,
	and the expression 
	for the controllability time 
	is valid in higher dimension.
	The main technique used 
	in this paper consists 
	in considering auxiliary systems obtained via orthogonal projection.	 
     Our results are exact in the case when 
     the range of $B$ is of co-dimension one (Theorem \ref{thm dim H = 1}).      
     As shown in Example \ref{example bar T ne T}, 
     the conclusions of 
     Theorem \ref{thm dim H = 1} do not hold 
     without the assumption $\dim H = 1$.

	We conclude with the following remarks about 
	desirable extensions.

	\begin{itemize}
		\item In this paper we worked with 
		a convex bounded constraint set.
		% Naturally 
		%one would want to know about 
		%a more general constraint set. We start by saying that some 
		%minor extensions seem to be within the reach
		%of the ideas and techniques developed in this paper.
		%In particular,
		In our analysis
		we can replace the assumption 
		that $\C$ is convex
		 with the assumption that 
		the projection of $\C$ on $H$
		is convex.
		On the other hand, unbounded 
		constraint sets require further considerations.
		We also  
		expect
		the ideas developed here
		to be applicable 
		to the case of a nice bounded constraint set with holes,
		 for example an annulus.
		
		\item The case when $\dim H \geq 2$
		is intriguing. Lemma \ref{T_C leq bar T_C} and Proposition \ref{T leq tilde T} 
		do shed some light on relation between  equations \eqref{the model}
		and \eqref{system hat}, whereas 
		Example \ref{example bar T ne T} demonstrates pitfalls 
		of trying to express $T_\C$ via $\overline T_\C$.
		Our intuitive guess is that Example \ref{example bar T ne T}
		is rather contrived, and `in most cases' 
		the equality $ T_\C = \overline T _\C$ 
		should hold. This `in most cases' 
		could be formulated as
		certain parameters of the system being generic
		as in \cite[Theorem 1]{LTZ18} (that is, belonging to 
		a dense open set), 
		or, alternatively, 
		as $ T_\C = \overline T _\C$  
		holding with probability one 
		when the parameters of the system
		are drawn from some continuous distributions.

	\end{itemize}

	\section*{Appendix}
	
	Here we formulate and prove a technical result
	used in Section \ref{sec theoretical part}. 
	The next lemma is used to establish a lower bound for the time $T$
	the system needs to travel a certain distance $M$.

	\begin{manualLem}{A1}\label{speedBound}
		Let $f$ be a non-negative  differentiable function 
		with $f(0) = 0$, $\dot f(t) \leq g(f(t))$, and $f(T)  = M$
		for $T, M > 0$ and a positive function $g$.
		Then $T \geq \int\limits _0 ^M \frac{dv}{g(v)}$,
		and the equality is achieved if $\dot f(t) = g(f(t))$.
	\end{manualLem}
	\textbf{Proof}.
	Let $h$ be defined as the solution to 
	\begin{equation}\label{stork}
	\dot h (t) = g(h(t)),  \  \ \  { h(0) = 0}.
	\end{equation}

	By the comparison theorem, $f(t) \leq h(t) $, $t \geq 0$.
	It follows from \eqref{stork} that $h$ is strictly increasing and thus inversible,
	and hence
	$\frac{dt}{d{  h}} = \frac{1}{g(h)} $, $t = \int\limits _0 ^h \frac{dv}{g(v)}$.
	Hence for $T_M:=   \int\limits _0 ^M \frac{dv}{g(v)}$ we get 
	$h(T_M) = M$. Thus,
	$f(T_M ) \leq h(T_M) = M$, and therefore
	$f(T) = M$ implies $T \geq T_M = \int\limits _0 ^M \frac{dv}{g(v)}$.
	Finally, if $ f = h$, then $T =  T_M = \int\limits _0 ^M \frac{dv}{g(v)}$.
	\qed
	
	\section*{Acknowledgements}

	This work has been partially supported by the project of the Italian Ministry of Education, Universities and Research (MIUR) ``Dipartimenti di Eccellenza 2018-2022".
	The authors thank the anonymous referees 
	for their insightful and careful reviews.
	
	%\section*{References}
	
	\bibliographystyle{alpha}
	\bibliography{SHS}

\end{document}